\documentclass[a4paper,oneside,11 pt]{article}%
\usepackage{amsmath}
\usepackage{amsfonts}
\usepackage{amssymb}
\usepackage{graphicx}%
\providecommand{\U}[1]{\protect\rule{.1in}{.1in}}

\newtheorem{theorem}{Theorem}[section]

\newtheorem{definition}[theorem]{Definition}

\newtheorem{lemma}[theorem]{Lemma}

\numberwithin{equation}{section}

\begin{document}

\title{The Application of G-heat equation and Numerical Properties \thanks{This work
was supported by National Natural Science Foundation of China (No. 11171187,
No. 10871118 and No. 10921101).}}
\author{Xiaolin Gong \thanks{Institute for Economics and Institute of Financial
Studies, Shandong University, Jinan, Shandong 250100, PR
China(agcaelyn@gmail.com).}
\and Shuzhen Yang \thanks{School of mathematics, Shandong University, Jinan,
Shandong 250100, PR China(yangsz@mail.sdu.edu.cn).}}
\date{}
\maketitle

\textbf{Abstract}: We consider a nonlinear expectation G-expectation which was
established by Peng. In order to compute the nonlinear probability under the
G-expectation, we prove that a function (special point is the nonlinear
probability) is the viscosity solution of the G-heat equation, and\ show that
the fully implicit discretization convergence to the viscosity solution of the
G-heat equation.


\addcontentsline{toc}{section}{\hspace*{1.8em}Abstract}

\section{Introduction}

In the mathematical Finance, we focus on the compute of probability of
default. Under the assumption of\ linear probability (expectation) space, we
use log normal distribution to describe the return of stock, and we could
easily\ calculus probability of default by normal distribution. For general
case, there is not only one probability. We need introduce volatility
uncertainty (including much more probabilities) in the market.

A nonlinear expectation (probability) G-expectation was established by Peng in
recent years, which could\ be equivalent to a set of probabilities (see
\cite{DHP11}). In the theory of G-expectation, the G-normal distribution and
G-Brownian motion were introduced and the corresponding stochastic calculus of
Ito's type were established (see \cite{P07a}, \cite{P08a}, \cite{P09}). In
Markovian case, the G-expectation is associated with fully nonlinear PDEs, and
is applied among economic and financial models with volatility uncertainty
(see \cite{Pooey}).

The next equation is used to compute the nonlinear probability (\cite{P07a}):%

\begin{equation}
\partial_{t}u-\frac{1}{2}(\bar{\sigma}^{2}(D_{xx}^{2}u)^{+}-\underline{\sigma
}^{2}(D_{xx}^{2}u)^{-})=0,\ u(0,x)=\varphi(x),x\in R, \label{in-1}%
\end{equation}

where $\varphi(x)=1_{\{x<0\}},x\in R.$

We show that $u(t,x):=\mathbb{\hat{E}}[\varphi(x+\sqrt{t}X)]$, $(t,x)\in
\lbrack0,\infty)\times R^{d}$, is the viscosity\ solution of the equation
(\ref{in-1}), where $\mathbb{\hat{E}}$ is the nonlinear expectation.

Following the work of\ \cite{Pooey}, \cite{Pang}, \cite{Qi}, \cite{Sun}, we
prove that the the fully implicit discretization convergence to the viscosity
solution of the G-heat equation.

Under the same maximum volatility, we compare the nonlinear probability
$u(1,0)$ and linear probability $\hat{u}(1,0)$ of the next two equations:
\begin{equation}%
\begin{array}
[c]{l}%
\partial_{t}u-\frac{1}{2}((D_{xx}^{2}u)^{+}-\frac{1}{4}(D_{xx}^{2}u)^{-})=0,\\
u(0,x)=I_{x\leq0},\text{ \ \ }x\in\lbrack-10,10],\\
u(t,-10)=1,\text{ }u(t,10)=0,\text{\ \ }t\in\lbrack0,1];
\end{array}
\label{in-2}%
\end{equation}

and%

\begin{equation}%
\begin{array}
[c]{l}%
\partial_{t}\hat{u}-\frac{1}{2}((D_{xx}^{2}\hat{u})^{+}-(D_{xx}^{2}\hat
{u})^{-})=0,\\
\hat{u}(0,x)=I_{x\leq0},\text{ \ \ }x\in\lbrack-10,10],\\
\hat{u}(t,-10)=1,\text{ }\hat{u}(t,10)=0,\text{\ \ }t\in\lbrack0,1].
\end{array}
\label{in-3}%
\end{equation}

By calculation, we have
\[
\mathbb{\hat{E}}[I_{X\leq0}]=u(1,0)=0.6680,\text{ }P(X\leq0)=\hat
{u}(1,0)=0.5010.\text{\ }%
\]

The paper is organized as follows: In section 2, the notations and results on
G-expectation is presented. In section 3, we prove that a funtion is the
viscosity solution of the G-heat equation. Then, we compare the value of
nonlinear probability and linear probability\ in section 4. The fully
implicit\ numerical convergence to the viscosity solution of the G-heat
equation is established in section 5.

\section{Preliminaries}

Firstly, we give the basic theory of G-expectation.

Let $\Omega$ be a given set and $\mathcal{H}$ a vector lattice of real valued
functions defined on $\Omega$, namely $c\in\mathcal{H}$ for each constant $c$
and $|X|\in\mathcal{H}$ if $X\in\mathcal{H}$. $\mathcal{H}$ is considered as
the space of random variables.

\begin{definition}
\label{def2.1} A sublinear expectation $\mathbb{\hat{E}}$ on $\mathcal{H}$ is
a functional $\mathbb{\hat{E}}:\mathcal{H}\rightarrow R$ satisfying the
following properties: for all $X,Y\in\mathcal{H}$, we have

(a) Monotonicity: If $X\geq Y$ then $\mathbb{\hat{E}}[X]\geq\mathbb{\hat{E}%
}[Y]$;

(b) Constant preservation: $\mathbb{\hat{E}}[c]=c$;

(c) Sub-additivity: $\mathbb{\hat{E}}[X+Y]\leq\mathbb{\hat{E}}[X]+\mathbb{\hat
{E}}[Y]$;

(d) Positive homogeneity: $\mathbb{\hat{E}}[\lambda X]=\lambda\mathbb{\hat{E}%
}[X]$ for each $\lambda\geq0$.

\begin{description}
\item $(\Omega,\mathcal{H},\mathbb{\hat{E}})$ is called a sublinear
expectation space.
\end{description}
\end{definition}

\begin{definition}
\label{def2.2} Let $X_{1}$ and $X_{2}$ be two $n$-dimensional random vectors
defined respectively in sublinear expectation spaces $(\Omega_{1}%
,\mathcal{H}_{1},\mathbb{\hat{E}}_{1})$ and $(\Omega_{2},\mathcal{H}%
_{2},\mathbb{\hat{E}}_{2})$. They are called identically distributed, denoted
by $X_{1}\overset{d}{=}X_{2}$, if $\mathbb{\hat{E}}_{1}[\varphi(X_{1}%
)]=\mathbb{\hat{E}}_{2}[\varphi(X_{2})]$, for all$\ \varphi\in C_{l.Lip}%
(R^{n})$, where $C_{l.Lip}(R^{n})$ be the space of real continuous functions
defined on $R^{n}$ such that
\[
|\varphi(x)-\varphi(y)|\leq C(1+|x|^{k}+|y|^{k})|x-y|\ \text{\ for
all}\ x,y\in R^{n},
\]
where $k$ and $C$ depend only on $\varphi$.
\end{definition}

\begin{definition}
\label{def2.3} In a sublinear expectation space $(\Omega,\mathcal{H}%
,\mathbb{\hat{E}})$, a random vector $Y=(Y_{1},\cdots,Y_{n})$, $Y_{i}%
\in\mathcal{H}$, is said to be independent of another random vector
$X=(X_{1},\cdots,X_{m})$, $X_{i}\in\mathcal{H}$ under $\mathbb{\hat{E}}%
[\cdot]$, denoted by $Y\bot X$, if for every test function $\varphi\in
C_{l.Lip}(R^{m}\times R^{n})$ we have $\mathbb{\hat{E}}[\varphi
(X,Y)]=\mathbb{\hat{E}}[\mathbb{\hat{E}}[\varphi(x,Y)]_{x=X}]$.
\end{definition}

\begin{definition}
\label{def2.4} ($G$-normal distribution) A $d$-dimensional random vector
$X=(X_{1},\cdots,X_{d})$ in a sublinear expectation space $(\Omega
,\mathcal{H},\mathbb{\hat{E}})$ is called $G$-normal distributed if for each
$a,b\geq0$ we have
\[
aX+b\bar{X}\overset{d}{=}\sqrt{a^{2}+b^{2}}X,
\]
where $\bar{X}$ is an independent copy of $X$, i.e., $\bar{X}\overset{d}{=}X$
and $\bar{X}\bot X$. Here the letter $G$ denotes the function
\[
G(A):=\frac{1}{2}\mathbb{\hat{E}}[\langle AX,X\rangle]:\mathbb{S}%
_{d}\rightarrow R,
\]
where $\mathbb{S}_{d}$ denotes the collection of $d\times d$ symmetric matrices.
\end{definition}

Peng \cite{P08a} showed that $X=(X_{1},\cdot\cdot\cdot,X_{d})$ is $G$-normally
distributed if and only if for each $\varphi\in C_{l.Lip}(R^{d})$,
$u(t,x):=\mathbb{\hat{E}}[\varphi(x+\sqrt{t}X)]$, $(t,x)\in\lbrack
0,\infty)\times R^{d}$, is the viscosity\ solution of the following $G$-heat
equation:%
\begin{equation}
\partial_{t}u-G(D_{xx}^{2}u)=0,\ u(0,x)=\varphi(x). \label{equa-1}%
\end{equation}

The function $G(\cdot):\mathbb{S}_{d}\rightarrow R$ is a monotonic, sublinear
mapping on $\mathbb{S}_{d}$ and $G(A)=\frac{1}{2}\mathbb{\hat{E}}%
[(AX,X)]\leq\frac{1}{2}|A|\mathbb{\hat{E}}[|X|^{2}]=:\frac{1}{2}|A|\bar
{\sigma}^{2}$ implies that there exists a bounded, convex and closed subset
$\Gamma\subset\mathbb{S}_{d}^{+}$ such that
\[
G(A)=\frac{1}{2}\sup_{\gamma\in\Gamma}\text{\textrm{tr}}[\gamma A],
\]
where $\mathbb{S}_{d}^{+}$ denotes the collection of nonnegative elements in
$\mathbb{S}_{d}$.

Let $\{W_{t}\}$ be a classical $d$-dimensional Brownian motion on a
probability space $(\Omega^{0},\mathcal{F}^{0},P^{0})$ and let$\ F^{0}%
=\{\mathcal{F}_{t}^{0}\}$ be the augmented filtration generated by $W$. Set
\[
\mathcal{P}_{M}:=\{P_{\theta}:P_{\theta}=P^{0}\circ(B_{.}^{\theta,0}%
)^{-1},B_{t}^{\theta,0}=\int_{0}^{t}\theta_{s}dW_{s},\theta\in L_{F^{0}}%
^{2}([0,T];\Gamma)\},
\]
where $L_{F^{0}}^{2}([0,T];\Gamma)$ is the collection of $F^{0}$-adapted
square integrable measurable processes with values in $\Gamma$. Set
$\mathcal{P=}\overline{\mathcal{P}_{M}}\mathcal{\ }$the closure of
$\mathcal{P}_{M}$ under the topology of weak convergence, then $\mathcal{P}$
is weakly compact. \cite{DHP11} proved that $\mathcal{P}$ represents
$\mathbb{\hat{E}}$ on $L_{G}^{1}(\Omega_{T})$.

Let $d=1$, we\ consider the finite difference method to the next G-heat equation:%

\begin{equation}%
\begin{array}
[c]{l}%
\partial_{t}u-\frac{1}{2}(\bar{\sigma}^{2}(D_{xx}^{2}u)^{+}-\underline{\sigma
}^{2}(D_{xx}^{2}u)^{-})=0,\text{ \ \ }x\in R,\text{ \ }t>0,\\
u(0,x)=\varphi(x),\text{ \ \ }x\in R,
\end{array}
\label{equa-2}%
\end{equation}
which%
\[
(D_{xx}^{2}u)^{+}=%
\genfrac{\{}{.}{0pt}{}{D_{xx}^{2}u,\text{ \ \ }D_{xx}^{2}u\geq
0,}{0,\text{\ \ \ \ \ \ \ }D_{xx}^{2}u<0,}%
\]
and%
\[
(D_{xx}^{2}u)^{-}=%
\genfrac{\{}{.}{0pt}{}{-D_{xx}^{2}u,\text{ \ \ }D_{xx}^{2}u\leq
0,}{0,\text{\ \ \ \ \ \ \ \ }D_{xx}^{2}u>0.}%
\]

For generally, we focus on the case $\varphi(x)=1_{\{x<y\}},y\in R.$ Set
$\bar{\sigma}=\underline{\sigma}=\sigma_{0},$ then $u(1,0)=P(X<y)$,
$X\overset{d}{=}N(0,\sigma_{0}^{2})$,\ specially. Next, we consider the
function: $\varphi(x)=1_{\{x<0\}},x\in R.$

\section{The viscosity solution of G-heat equation}

We will show that $u(t,x):=\mathbb{\hat{E}}[\varphi(x+\sqrt{t}X)]$,
$(t,x)\in\lbrack0,\infty)\times R^{d}$, is the viscosity\ solution of the
following $G$-heat equation:%
\begin{equation}
\partial_{t}u-G(D_{xx}^{2}u)=0,\ u(0,x)=\varphi(x),x\in R, \label{vis-1}%
\end{equation}

where $\varphi(x)=1_{\{x<0\}},x\in R.$

\begin{lemma}
\label{lem1} $\lim_{n\rightarrow\infty}\mathbb{\hat{E}}[\varphi_{n}%
(X)]=\mathbb{\hat{E}}[\varphi(X)],$
\end{lemma}

where
\[
\varphi_{n}(x)=\left\{
\begin{array}
[c]{c}%
1,\text{ \ \ \ \ \ \ \ \ \ }x\leq0\\
1-nx,0<x<\frac{1}{n}\\
0,\text{ \ \ \ \ \ \ \ }\frac{1}{n}\leq x
\end{array}
.\right.
\]
\textbf{Proof. }For $\varphi(x)=1_{\{x<0\}},$ $\mathbb{\hat{E}}[\varphi
_{n}(X)]-\mathbb{\hat{E}}[\varphi(X)]\leq\mathbb{\hat{E}[}1_{\{0<X<\frac{1}%
{n}\}}\mathbb{]},$ and $X$ is a $G$-normal distribution, we have
\[%
\begin{array}
[c]{rl}%
\mathbb{\hat{E}[}1_{\{0<X<\frac{1}{n}\}}\mathbb{]} & =\sup_{p_{\theta}%
\in\mathcal{P}_{M}}P_{\theta}\mathbb{(}0<X<\frac{1}{n}\mathbb{)}\\
& =\sup_{\theta\in L_{F^{0}}^{2}([0,T];\Gamma)}P^{0}\mathbb{(}0<\int_{0}%
^{1}\theta_{s}dW_{s}<\frac{1}{n}\mathbb{)}\\
& \leq P^{0}\mathbb{(-}\frac{\bar{\sigma}}{\underline{\sigma}}\frac{1}{n}%
<\bar{\sigma}W_{1}<\frac{\bar{\sigma}}{\underline{\sigma}}\frac{1}%
{n}\mathbb{)}.
\end{array}
\]

So $\mathbb{\hat{E}}[\varphi_{n}(X)]-\mathbb{\hat{E}}[\varphi
(X)]\longrightarrow0,$ as $n\longrightarrow\infty.$

\begin{theorem}
\label{w11}The function $u(t,x):=\mathbb{\hat{E}}[\varphi(x+\sqrt{t}X)]$,
$(t,x)\in\lbrack0,\infty)\times R^{d}$, is the viscosity solution of equation
(\ref{vis-1}).
\end{theorem}

\textbf{Proof. }Firstly, we show that $u$ is continuous in $[0,\infty)\times
R.$

For $\forall\delta>0,$ $u(t+\delta,x)-u(t,x)=\mathbb{\hat{E}}[\varphi
(x+\sqrt{t+\delta}X)]-\mathbb{\hat{E}}[\varphi(x+\sqrt{t}X)],$ by Lemma
\ref{lem1} and the definition of $\mathcal{P}_{M}:$%
\[%
\begin{array}
[c]{rl}
& \mathbb{\hat{E}}[\varphi(x+\sqrt{t+\delta}X)]-\mathbb{\hat{E}}%
[\varphi(x+\sqrt{t}X)]\\
= & \lim_{n\rightarrow\infty}(\mathbb{\hat{E}}[\varphi_{n}(x+\sqrt{t+\delta
}X)]-\mathbb{\hat{E}}[\varphi_{n}(x+\sqrt{t}X)])\\
= & \lim_{n\rightarrow\infty}(\mathbb{\hat{E}}[\varphi_{n}(x+\sqrt{\delta}%
\bar{X}+\sqrt{t}X)]-\mathbb{\hat{E}}[\varphi_{n}(x+\sqrt{t}X)])\\
= & \mathbb{\hat{E}}[\varphi(x+\sqrt{\delta}\bar{X}+\sqrt{t}X)]-\mathbb{\hat
{E}}[\varphi(x+\sqrt{t}X)]\\
\leq & \sup_{p_{\theta}\in\mathcal{P}_{M}}P_{\theta}(\{x+\sqrt{t}X\leq
-\sqrt{\delta}\bar{X}\}/\{x+\sqrt{t}X\leq0\})\\
\leq & \sup_{\theta\in L_{F^{0}}^{2}([0,T];\Gamma)}[P^{0}(0<x+\sqrt{t}\int%
_{0}^{1}\theta_{s}dW_{s}\leq-\sqrt{\delta}\int_{0}^{1}\theta_{s}d\bar{W}%
_{s})\\
& +P^{0}(-\sqrt{\delta}\int_{0}^{1}\theta_{s}d\bar{W}_{s}<x+\sqrt{t}\int%
_{0}^{1}\theta_{s}dW_{s}\leq0)],
\end{array}
\]

where $\bar{W}$ is independent identically distributed with $W$ in the linear
expectation space $(\Omega^{0},\mathcal{F}^{0},P^{0}).$

By simple calculus, we have%

\[%
\begin{array}
[c]{cl}
& \sup_{\theta\in L_{F^{0}}^{2}([0,T];\Gamma)}P^{0}(0<x+\sqrt{t}\int_{0}%
^{1}\theta_{s}dW_{s}\leq-\sqrt{\delta}\int_{0}^{1}\theta_{s}d\bar{W}_{s})\\
\leq & P^{0}(\frac{\bar{\sigma}}{\underline{\sigma}}\sqrt{\delta}\bar{\sigma
}\bar{W}_{1}<\frac{\underline{\sigma}}{\bar{\sigma}}x+\frac{\underline{\sigma
}}{\bar{\sigma}}\sqrt{t}\bar{\sigma}W_{1}\leq-\frac{\bar{\sigma}%
}{\underline{\sigma}}\sqrt{\delta}\bar{\sigma}\bar{W}_{1}),
\end{array}
\]

then, by dominated convergence theorem,
\[
\sup_{\theta\in L_{F^{0}}^{2}([0,T];\Gamma)}P^{0}(0<x+\sqrt{t}\int_{0}%
^{1}\theta_{s}dW_{s}\leq-\sqrt{\delta}\int_{0}^{1}\theta_{s}d\bar{W}%
_{s})\longrightarrow0,\text{ as }\delta\longrightarrow0.
\]
So $u$ is continuous in $t.$ Similarly we could prove that $u$ is continuous
in $x.$

Now, we prove that $u$ is a viscosity subsolution of equation\ (\ref{vis-1}). \

For a fixed $(t,x)\in$ $[0,\infty)\times R^{d},$ let $\psi\in C_{b}%
^{2,3}([0,\infty)\times R^{d}),$ such that $\psi\geq u$ and $\psi
(t,x)=u(t,x).$ By Taylor's expansion, it follows that, for $\delta\in(0,t),$%

\[%
\begin{array}
[c]{cl}%
0 & \leq\mathbb{\hat{E}[\psi(}t-\mathbb{\delta},x+\sqrt{\delta}%
X\mathbb{)-\mathbb{\psi(}}t,x\mathbb{\mathbb{)}]}\\
& \leq o(\delta)-\partial_{t}\psi(t,x)\delta+\mathbb{\hat{E}}[\langle
D_{x}\psi(t,x),X\rangle\sqrt{\delta}+\frac{1}{2}\langle D_{xx}^{2}%
\psi(t,x)X,X\rangle\delta]\\
& =o(\delta)-\partial_{t}\psi(t,x)\delta+\mathbb{\hat{E}}[\frac{1}{2}\langle
D_{xx}^{2}\psi(t,x)X,X\rangle\delta]\\
& =-\partial_{t}\psi(t,x)\delta+\delta G(D_{xx}^{2}\psi)(t,x)+o(\delta),
\end{array}
\]

so
\[
\partial_{t}\psi(t,x)-G(D_{xx}^{2}\psi)(t,x)\leq0.
\]

Thus $u$ is a viscosity subsolution of (\ref{vis-1}). Similarly, we have $u$
is a viscosity solution of (\ref{vis-1}).\

This completes the proof.

\section{Numerical Example}

In this section, we give an example which is important for financial market.
By Theorem \ref{w11}, the nonlinear probability $u(t,x):=\mathbb{\hat{E}%
}[I_{x+\sqrt{t}X\leq0}]$, $(t,x)\in\lbrack0,\infty)\times R$, is the viscosity
solution of the following $G$-heat equation:%
\begin{equation}
\partial_{t}u-G(D_{xx}^{2}u)=0,\ u(0,x)=I_{x\leq0}. \label{Nequa-1}%
\end{equation}

i.e.,%

\begin{equation}%
\begin{array}
[c]{l}%
\partial_{t}u-\frac{1}{2}(\bar{\sigma}^{2}(D_{xx}^{2}u)^{+}-\underline{\sigma
}^{2}(D_{xx}^{2}u)^{-})=0,\text{ \ \ }x\in R,\text{ \ }t>0,\\
u(0,x)=I_{x\leq0},\text{ \ \ }x\in R.
\end{array}
\label{Nequa-2}%
\end{equation}

Next, we consider a boundary\ problem of (\ref{Nequa-2}), i.e.,
\begin{equation}%
\begin{array}
[c]{l}%
\partial_{t}u-\frac{1}{2}((D_{xx}^{2}u)^{+}-\frac{1}{4}(D_{xx}^{2}u)^{-})=0,\\
u(0,x)=I_{x\leq0},\text{ \ \ }x\in\lbrack-10,10],\\
u(t,-10)=1,\text{ }u(t,10)=0,\text{\ \ }t\in\lbrack0,1].
\end{array}
\label{Nequa-3}%
\end{equation}

For a given probability space $(\Omega,\mathcal{F},P),$ the linear probability
$u(t,x):=E[I_{x+\sqrt{t}X\leq0}]$, $(t,x)\in\lbrack0,\infty)\times R$, is the
viscosity solution of the following heat equation:%

\begin{equation}
\partial_{t}u-D_{xx}^{2}u=0,\ u(0,x)=I_{x\leq0}. \label{Nequa-4}%
\end{equation}

We also\ consider a boundary\ problem of (\ref{Nequa-4}), i.e.,%

\begin{equation}%
\begin{array}
[c]{l}%
\partial_{t}\hat{u}-\frac{1}{2}((D_{xx}^{2}\hat{u})^{+}-(D_{xx}^{2}\hat
{u})^{-})=0,\\
\hat{u}(0,x)=I_{x\leq0},\text{ \ \ }x\in\lbrack-10,10],\\
\hat{u}(t,-10)=1,\text{ }\hat{u}(t,10)=0,\text{\ \ }t\in\lbrack0,1].
\end{array}
\label{Nequa-5}%
\end{equation}

Comparing the value of $u(1,x)$ and $\hat{u}(1,x),x\in\lbrack-10,10]:$

\begin{center}
\includegraphics[width=3.3 in]{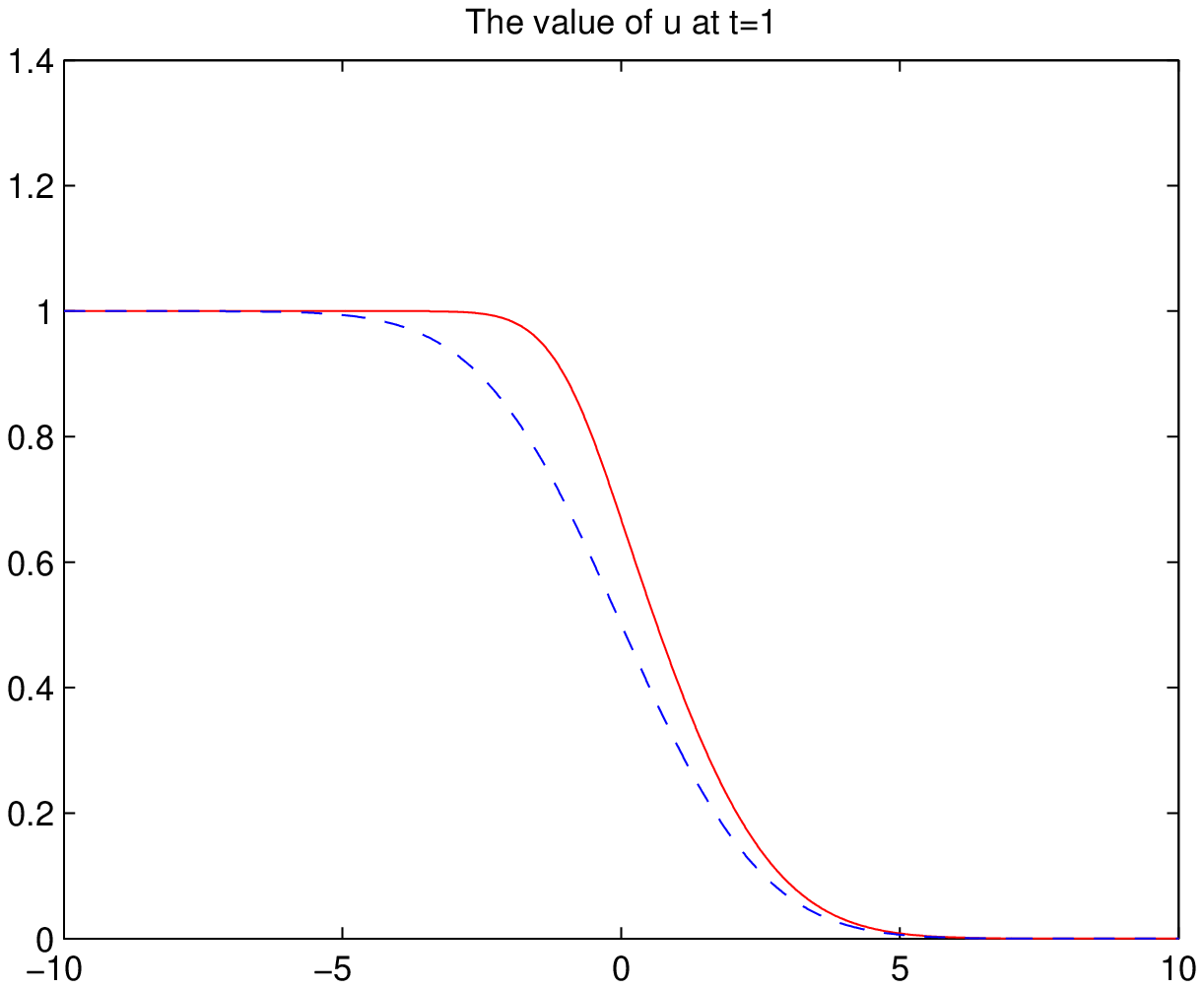}
\end{center}

The red line is the value of function $u(1,x),$ $x\in\lbrack-10,10]$, and the
blue line is the the value of function $\hat{u}(1,x),$ $x\in\lbrack-10,10]$.
By $u(t,x):=\mathbb{\hat{E}}[I_{x+\sqrt{t}X\leq0}]$ and $\hat{u}%
(1,x)=E[I_{x+\sqrt{t}X\leq0}],$ we have%

\[
\mathbb{\hat{E}}[I_{X\leq0}]=u(1,0)=0.6680,\text{ }P(X\leq0)=\hat
{u}(1,0)=0.5010.\text{\ }%
\]

\section{Numerical Analysis}

In this section, we consider the bounded boundary\ problem of (\ref{equa-2}),
i.e.,
\begin{equation}%
\begin{array}
[c]{l}%
\partial_{t}u-\frac{1}{2}(\bar{\sigma}^{2}(D_{xx}^{2}u)^{+}-\underline{\sigma
}^{2}(D_{xx}^{2}u)^{-})=0,\\
u(0,x)=\varphi(x),\text{ \ \ }x\in\lbrack a,b],\\
u(t,a)=g(t),\text{ }u(t,b)=h(t)\text{\ \ }t\in\lbrack0,T].
\end{array}
\label{equa-3}%
\end{equation}

where $\varphi,g,h$ are bounded and measurable functions.

\subsection{A Finite Difference Discretization}

The equation (\ref{equa-3}) can be discretized by a standard finite difference
method with variable timeweighting to give%

\begin{equation}%
\begin{array}
[c]{cl}%
u_{i}^{n+1}-u_{i}^{n}= & \theta\alpha_{i}^{n}[u_{i+1}^{n}-2u_{i}^{n}%
+u_{i-1}^{n}]+(1-\theta)\alpha_{i}^{n+1}[u_{i+1}^{n+1}-2u_{i}^{n+1}%
+u_{i-1}^{n+1}],
\end{array}
\label{FDD-1}%
\end{equation}

where
\begin{equation}%
\begin{array}
[c]{rl}%
\alpha_{i}^{n}:= & \frac{\sigma(\Gamma_{i}^{n})^{2}\triangle t_{i}}%
{2(x_{i+1}-x_{i})(x_{i}-x_{i-1})},\\
\sigma(\Gamma_{i}^{n}):= & \left\{
\begin{array}
[c]{l}%
\bar{\sigma},\text{ \ \ if }\Gamma_{i}^{n}\geq0\\
\underline{\sigma},\text{ \ \ if }\Gamma_{i}^{n}<0
\end{array}
\right.  ,\\
\Gamma_{i}^{n}:= & \frac{u_{i+1}^{n}-2u_{i}^{n}+u_{i-1}^{n}}{(x_{i+1}%
-x_{i})(x_{i}-x_{i-1})}.
\end{array}
\label{FDD-2}%
\end{equation}

In this paper, we consider the fully implicit schemes\ with $\theta=0,$ i.e.,%

\begin{equation}%
\begin{array}
[c]{cl}%
u_{i}^{n+1}-u_{i}^{n}= & \alpha_{i}^{n+1}[u_{i+1}^{n+1}-2u_{i}^{n+1}%
+u_{i-1}^{n+1}],
\end{array}
\label{PDD-3}%
\end{equation}

The set of algebraic equation (\ref{equa-3}) is nonlinear for the formula of
$\sigma(\Gamma_{i}^{n}).$ So we consider the discrete equation at each node as%

\[%
\begin{array}
[c]{c}%
\varphi_{i}^{n}:=u_{i}^{n}-u_{i}^{n+1}+\alpha_{i+1}^{n+1}[u_{i+1}^{n+1}%
-2u_{i}^{n+1}+u_{i-1}^{n+1}].
\end{array}
\]

Following the work of D.M. Pooey \cite{Pooey} (more see, Pang and Qi
\cite{Pang}; Qi and Sun \cite{Qi}; Sun and Han \cite{Sun}), we must specify
the element of the generalized Jacobian that will be used in the Newton
iteration. We define the derivatives as%

\[
\frac{\partial\sigma(\Gamma)^{2}\Gamma}{\partial\Gamma}=\left\{
\begin{array}
[c]{l}%
\bar{\sigma}^{2},\text{ \ \ if }\Gamma\geq0\\
\underline{\sigma}^{2},\text{ \ \ if }\Gamma<0
\end{array}
\right.  ,
\]

For further analysis the Newton iteration, we rewrite the discrete equation
(\ref{FDD-1}) in matrix form. Let%

\[%
\begin{array}
[c]{l}%
U^{n+1}=[u_{0}^{n+1},u_{1}^{n+1},\cdots,u_{m}^{n+1}]^{\prime},\text{ \ }%
U^{n}=[u_{0}^{n},u_{1}^{n},\cdots,u_{m}^{n}]^{\prime},\\
\lbrack M^{n}U^{n}]_{i}:=-\alpha_{i}^{n}[u_{i+1}^{n}-2u_{i}^{n}+u_{i-1}^{n}].
\end{array}
\]

For convenience, we modify the first and last rows of $M$ as needed to handle
the bounded boundary conditions. By the discretization schemes in
(\ref{PDD-3}), the matrix $M$ is a diagonally dominant matric with positive
diagonals and non-positive off-diagonals. Note that all the elements of the
inverse of $M$ are non-negative. The discrete equation (\ref{PDD-3}) can be
rewritten as:%

\begin{equation}
\lbrack I+M^{n+1}]U^{n+1}=U^{n}, \label{PDD-4}%
\end{equation}

where $I$ is the identity matrix. Next we prove the convergence of the\ Newton
iteration for full implicit schemes.

\subsection{Convergence of the Newtion Iteration\ Schemes}

For the matrix $M$ is a diagonally dominant matric, we can analysis the Newton
iteration of equation (\ref{PDD-4}). We adopt the Newton timestep as the
following scheme:

(a) Let $(U^{n+1})^{0}=U^{n};$

(b) For $k=0,1,2,\cdots$ Solve%

\begin{equation}
\lbrack I+M((U^{n+1})^{k})](U^{n+1})^{k+1}=U^{n} \label{PDD-5}%
\end{equation}
where $(U^{n+1})^{k+1}$ is the $(k+1)$th iteration, and $M((U^{n+1})^{k})$
means $M$ be dependent on $(U^{n+1})^{k}$.

(c) For a given small number $\varepsilon$, if
\[
\max_{i}\left\vert (u_{i}^{n+1})^{k+1}-(u_{i}^{n+1})^{k}\right\vert
<\varepsilon\cdot\max_{i}(1,\left\vert (u_{i}^{n+1})^{k+1}\right\vert ).
\]

(d) we end the scheme.

We show the convergence results about the above Newton iteration as follows:

\begin{theorem}
\label{w3}The nonlinear iteration (\ref{PDD-5}) convergence to the unique
solution of (\ref{PDD-4}), for given intial iterate $(U^{n+1})^{0}=U^{n}.$
\end{theorem}

\noindent\textbf{Proof.\ }For notional convergence, we denote $\hat{M}%
^{k}=M((U^{n+1})^{k})$ and $\hat{U}^{k}=(U^{n+1})^{k}.$ So
equation\ (\ref{PDD-5}) can be rewritten as%
\begin{equation}
\lbrack I+\hat{M}^{k}]\hat{U}^{k+1}=U^{n}. \label{PDD-6}%
\end{equation}

Firstly, we show that the sequence $\{\hat{U}^{k}\}_{0\leq k}$ is
monotonically. The $k$ iteration of equation (\ref{PDD-6}) gives
\begin{equation}
\lbrack I+\hat{M}^{k-1}]\hat{U}^{k}=U^{n}. \label{PDD-7}%
\end{equation}

Subtracting equation (\ref{PDD-6}) from equation (\ref{PDD-7}), we have%

\begin{equation}
\lbrack I+\hat{M}^{k}](\hat{U}^{k+1}-\hat{U}^{k})=[\hat{M}^{k-1}-\hat{M}%
^{k}]\hat{U}^{k}. \label{PDD-8}%
\end{equation}

We consider the right side of \ref{PDD-8}, for each $i$%

\[
([\hat{M}^{k-1}-\hat{M}^{k}]\hat{U}^{k})_{i}=\frac{\triangle t_{i}(\sigma
(\hat{\Gamma}_{i}^{k})^{2}-\sigma(\hat{\Gamma}_{i}^{k-1})^{2})}{2}\hat{\Gamma
}_{i}^{k},
\]

where%
\[%
\begin{array}
[c]{l}%
\hat{\Gamma}_{i}^{k}=\frac{\hat{u}_{i+1}^{k}-2\hat{u}_{i}^{k}+\hat{u}%
_{i-1}^{k}}{(x_{i+1}-x_{i})(x_{i}-x_{i-1})},\\
\hat{U}^{k}:=[\hat{u}_{0}^{k},\hat{u}_{1}^{k},\cdots,\hat{u}_{m}^{k}]^{\prime
},\\
\sigma(\hat{\Gamma}_{i}^{k}):=\left\{
\begin{array}
[c]{l}%
\bar{\sigma},\text{ \ \ if }\hat{\Gamma}_{i}^{k}\geq0\\
\underline{\sigma},\text{ \ \ if }\hat{\Gamma}_{i}^{k}<0
\end{array}
\right.  .
\end{array}
\]

By the equation (\ref{FDD-2}), if $\hat{\Gamma}_{i}^{k}\leq0,$ $\sigma
(\hat{\Gamma}_{i}^{k})^{2}=\underline{\sigma},$ then
\[
\frac{\triangle t_{i}(\sigma(\hat{\Gamma}_{i}^{k})^{2}-\sigma(\hat{\Gamma}%
_{i}^{k-1})^{2})}{2}\hat{\Gamma}_{i}^{k}\geq0;
\]

Similarly, if $\hat{\Gamma}_{i}^{k}\geq0,$ $\sigma(\hat{\Gamma}_{i}^{k}%
)^{2}=\bar{\sigma},$ then
\[
\frac{\triangle t_{i}(\sigma(\hat{\Gamma}_{i}^{k})^{2}-\sigma(\hat{\Gamma}%
_{i}^{k-1})^{2})}{2}\hat{\Gamma}_{i}^{k}\geq0.
\]

For the matric $I+\hat{M}^{k}$ is a diagonally dominant matric, the inverse of
matric $I+\hat{M}^{k}$ is non-negative, we have%

\begin{equation}
\hat{U}^{k+1}-\hat{U}^{k}\geq0,\text{ \ \ }k\geq1. \label{PDD-9}%
\end{equation}

Next, we need to prove the sequence $\{\hat{U}^{k}\}_{0\leq k}$ is bounded.
Set $C_{\max}=\max_{i}u_{i}^{n},$ $C_{\min}=\min_{i}u_{i}^{n},$ $\hat{U}%
_{\max}=\max_{i}\hat{u}_{i}^{k},$ $\hat{U}_{\min}=\min_{i}\hat{u}_{i}^{k}$. By
the equation (\ref{PDD-7}), we have%
\begin{equation}
\hat{u}_{i}^{k}-\hat{\alpha}_{i}^{k-1}[\hat{u}_{i+1}^{k}-2\hat{u}_{i}^{k}%
+\hat{u}_{i-1}^{k}]=u_{i}^{n}, \label{PDD-91}%
\end{equation}

where%
\begin{equation}
\alpha_{i}^{k}:=\frac{\sigma(\hat{\Gamma}_{i}^{k})^{2}\triangle t_{i}%
}{2(x_{i+1}-x_{i})(x_{i}-x_{i-1})}. \label{PDD-10}%
\end{equation}

By the equation (\ref{PDD-91}), and $\hat{\alpha}_{i}^{k-1}\geq0,$ then%

\[
(1+2\hat{\alpha}_{i}^{k-1})\hat{u}_{i}^{k}\leq2\hat{\alpha}_{i}^{k-1}\hat
{U}_{\max}+C_{\max},
\]

and
\[
(1+2\hat{\alpha}_{i}^{k-1})\hat{u}_{i}^{k}\geq2\hat{\alpha}_{i}^{k-1}\hat
{U}_{\min}+C_{\min}.
\]

So%

\[
\hat{u}_{i}^{k}\leq\frac{2\hat{\alpha}_{i}^{k-1}}{1+2\hat{\alpha}_{i}^{k-1}%
}\hat{U}_{\max}+C_{\max},
\]

and%
\[
\hat{u}_{i}^{k}\geq\frac{2\hat{\alpha}_{i}^{k-1}}{1+2\hat{\alpha}_{i}^{k-1}%
}\hat{U}_{\min}+C_{\min}.
\]

Set $\max_{i}\frac{2\hat{\alpha}_{i}^{k-1}}{1+2\hat{\alpha}_{i}^{k-1}}=b_{1},$
$\max_{i}\frac{2\hat{\alpha}_{i}^{k-1}}{1+2\hat{\alpha}_{i}^{k-1}}=b_{2},$ and
$0\leq b_{1},b_{2}<1,$ then we have%

\[
\hat{U}_{\max}\leq\frac{1}{1-b_{1}}C_{\max},\text{ \ }\hat{U}_{\min}\geq
\frac{1}{1-b_{2}}C_{\min}.
\]

Now, we prove the uniqueness. Suppose there are two solutions to equation
(\ref{PDD-4}), $U_{1}^{n+1}$ and $U_{2}^{n+1}$, such that%

\[
\lbrack I+M_{1}]U_{1}^{n+1}=U^{n},\text{ \ \ }[I+M_{2}]U_{2}^{n+1}=U^{n}.
\]

Similar the proof of the monotonicity sequence $\{\hat{U}^{k}\}_{0\leq k},$ we
have%
\[
\lbrack I+M_{2}](U_{2}^{n+1}-U_{1}^{n+1})=[M_{1}-M_{2}]U_{1}^{n+1},
\]

and
\[
U_{2}^{n+1}-U_{1}^{n+1}\geq0
\]

By the equality of $U_{1}^{n+1}$ and $U_{2}^{n+1},$ we have $U_{2}^{n+1}%
-U_{1}^{n+1}=0.$

Thus, we complete the proof.

\subsection{The Convergence of Fully Implicit}

In the above section. we have proved the convergence of the Newton iteration
for the nonlinear equation (\ref{PDD-3}). Next, we would to prove the full
implicit schemes convergence to the viscosity solution of (\ref{equa-3}). By
the work of Barles in \cite{Barles}, we know that a stable, consistent, and
monotone discretization will convergence to the viscosity solution.

\begin{theorem}
\label{w2} The fully implicit discretization (\ref{PDD-4}) convergences to the
solution of the equation (\ref{equa-3}), as $\triangle t$, $\triangle
x\rightarrow0.$
\end{theorem}

We first give some important lemmas for prove Theorem \ref{w2}.

Review the discrete equation at each node as%

\begin{equation}%
\begin{array}
[c]{c}%
\varphi_{i}^{n}:=u_{i}^{n}-u_{i}^{n+1}+\alpha_{i+1}^{n+1}[u_{i+1}^{n+1}%
-2u_{i}^{n+1}+u_{i-1}^{n+1}]
\end{array}
\label{CFI}%
\end{equation}

then at each step%

\begin{equation}
\varphi_{i}^{n}(u_{i+1}^{n+1},u_{i}^{n+1},u_{i-1}^{n+1},u_{i}^{n})=0,\text{
\ \ }\forall i. \label{CFI-1}%
\end{equation}

In the case fo nondifferentiable $\varphi_{i}^{n},$ we use the following
definition of monotonicity:

\begin{definition}
A discretization of the form (\ref{CFI-1}) is monotone if either%
\[%
\begin{array}
[c]{l}%
\varphi_{i}^{n}(u_{i+1}^{n+1}+\varepsilon_{i+1}^{n+1},u_{i}^{n+1}%
,u_{i-1}^{n+1}+\varepsilon_{i-1}^{n+1},u_{i}^{n}+\varepsilon_{i}^{n}%
)\geq\varphi_{i}^{n}(u_{i+1}^{n+1},u_{i}^{n+1},u_{i-1}^{n+1},u_{i}^{n}),\\
\varphi_{i}^{n}(u_{i+1}^{n+1},u_{i}^{n+1}+\varepsilon_{i}^{n+1},u_{i-1}%
^{n+1},u_{i}^{n})\leq\varphi_{i}^{n}(u_{i+1}^{n+1},u_{i}^{n+1},u_{i-1}%
^{n+1},u_{i}^{n}),\\
\forall\varepsilon_{i+1}^{n+1},\varepsilon_{i}^{n+1},\varepsilon_{i-1}%
^{n+1},\varepsilon_{i}^{n}\geq0,
\end{array}
\]

or%
\[%
\begin{array}
[c]{l}%
\varphi_{i}^{n}(u_{i+1}^{n+1}+\varepsilon_{i+1}^{n+1},u_{i}^{n+1}%
,u_{i-1}^{n+1}+\varepsilon_{i-1}^{n+1},u_{i}^{n}+\varepsilon_{i}^{n}%
)\leq\varphi_{i}^{n}(u_{i+1}^{n+1},u_{i}^{n+1},u_{i-1}^{n+1},u_{i}^{n}),\\
\varphi_{i}^{n}(u_{i+1}^{n+1},u_{i}^{n+1}+\varepsilon_{i}^{n+1},u_{i-1}%
^{n+1},u_{i}^{n})\geq\varphi_{i}^{n}(u_{i+1}^{n+1},u_{i}^{n+1},u_{i-1}%
^{n+1},u_{i}^{n}),\\
\forall\varepsilon_{i+1}^{n+1},\varepsilon_{i}^{n+1},\varepsilon_{i-1}%
^{n+1},\varepsilon_{i}^{n}\geq0.
\end{array}
\]

\end{definition}

Next, we prove the monotonicity of the fully implicit discretization.

\begin{lemma}
\label{w4}The fully implicit discretization (\ref{CFI}) is monotone, indepdent
of any choice of $\triangle t$ and $\triangle x.$
\end{lemma}

\noindent\textbf{Proof}: For any given $\varepsilon>0,$ we just to chek the
next two euation:%

\[%
\begin{array}
[c]{c}%
\varphi_{i}^{n}(u_{i+1}^{n+1}+\varepsilon,u_{i}^{n+1},u_{i-1}^{n+1},u_{i}%
^{n})\geq\varphi_{i}^{n}(u_{i+1}^{n+1},u_{i}^{n+1},u_{i-1}^{n+1},u_{i}^{n}),\\
\varphi_{i}^{n}(u_{i+1}^{n+1},u_{i}^{n+1}+\varepsilon,u_{i-1}^{n+1},u_{i}%
^{n})\leq\varphi_{i}^{n}(u_{i+1}^{n+1},u_{i}^{n+1},u_{i-1}^{n+1},u_{i}^{n}).
\end{array}
\]

By the definition of $\varphi_{i}^{n},$ we have%
\[%
\begin{array}
[c]{rl}%
\varphi_{i}^{n}(u_{i+1}^{n+1}+\varepsilon,u_{i}^{n+1},u_{i-1}^{n+1},u_{i}%
^{n})= & u_{i}^{n}-u_{i}^{n+1}+\alpha_{i+1}^{n+1}[u_{i+1}^{n+1}-2u_{i}%
^{n+1}+u_{i-1}^{n+1}]+\alpha_{i+1}^{n+1}\cdot\varepsilon\\
\geq & u_{i}^{n}-u_{i}^{n+1}+\alpha_{i+1}^{n+1}[u_{i+1}^{n+1}-2u_{i}%
^{n+1}+u_{i-1}^{n+1}]\\
= & \varphi_{i}^{n}(u_{i+1}^{n+1},u_{i}^{n+1},u_{i-1}^{n+1},u_{i}^{n})
\end{array}
\]

and%
\[%
\begin{array}
[c]{rl}%
\varphi_{i}^{n}(u_{i+1}^{n+1},u_{i}^{n+1}+\varepsilon,u_{i-1}^{n+1},u_{i}%
^{n})= & u_{i}^{n}-u_{i}^{n+1}+\alpha_{i+1}^{n+1}[u_{i+1}^{n+1}-2u_{i}%
^{n+1}+u_{i-1}^{n+1}]-(2\alpha_{i+1}^{n+1}+1)\cdot\varepsilon\\
\leq & u_{i}^{n}-u_{i}^{n+1}+\alpha_{i+1}^{n+1}[u_{i+1}^{n+1}-2u_{i}%
^{n+1}+u_{i-1}^{n+1}]\\
= & \varphi_{i}^{n}(u_{i+1}^{n+1},u_{i}^{n+1},u_{i-1}^{n+1},u_{i}^{n}).
\end{array}
\]

This completes the proof.

\noindent\textbf{Proof of Theorem \ref{w2}:}

By the results of Barles, we just to check that the fully
implicit\ discretization is consistent, stable, monotone. Fristly, the formula
(\ref{CFI-1}) is a consistent\textbf{\ }discretization. Then Theorem \ref{w3}
shows that the fully implicit\ discretization is monotone. So we need to prove
the\ discretization is stable. Set%

\[
U_{\max}^{n}=\max(\max_{i}U_{i}^{n}),g^{n},h^{n}),\text{ \ \ }U_{\min}%
^{n}=\min(\min_{i}U_{i}^{n}),g^{n},h^{n}).
\]
where $g^{n},h^{n}$ is the boundary value of the $n$th times step. Using the
same mathod as in Lemma \ref{w4}, we have the more exact results:%

\[
U_{\min}^{n}\leq u_{i}^{n+1}\leq U_{\max}^{n}.
\]

Thus, we complete the proof.

\subsection{The Superlinear Expectation}

For reader convenience, we still use the same notions as in sublinear
expectation (G-expectation), and show the main results of superlinear expectation.

\begin{definition}
\label{def2.1 copy(1)} A superlinear expectation $\mathbb{\hat{E}}$ on
$\mathcal{H}$ is a functional $\mathbb{\hat{E}}:\mathcal{H}\rightarrow R$
satisfying the following properties: for all $X,Y\in\mathcal{H}$, we have

(a) Monotonicity: If $X\geq Y$ then $\mathbb{\hat{E}}[X]\geq\mathbb{\hat{E}%
}[Y]$;

(b) Constant preservation: $\mathbb{\hat{E}}[c]=c$;

(c) Sub-additivity: $\mathbb{\hat{E}}[X+Y]\leq\mathbb{\hat{E}}[X]+\mathbb{\hat
{E}}[Y]$;

(d) Positive homogeneity: $\mathbb{\hat{E}}[\lambda X]=\lambda\mathbb{\hat{E}%
}[X]$ for each $\lambda\geq0$.

\begin{description}
\item $(\Omega,\mathcal{H},\mathbb{\hat{E}})$ is called a sublinear
expectation space.
\end{description}
\end{definition}

The bounded boundary\ problem is
\begin{equation}%
\begin{array}
[c]{l}%
\partial_{t}u-\frac{1}{2}(\underline{\sigma}^{2}(D_{xx}^{2}u)^{+}-\bar{\sigma
}^{2}(D_{xx}^{2}u)^{-})=0,\\
u(0,x)=\varphi(x),\text{ \ \ }x\in\lbrack a,b],\\
u(t,a)=g(t),\text{ }u(t,b)=h(t)\text{\ \ }t\in\lbrack0,T].
\end{array}
\label{sequa-3}%
\end{equation}

where $\varphi,g,h$ are measureable functions.

The equation (\ref{sequa-3}) can be discretized by a standard finite
difference method with variable timeweighting to give%

\begin{equation}%
\begin{array}
[c]{cl}%
u_{i}^{n+1}-u_{i}^{n}= & \theta\alpha_{i}^{n}[u_{i+1}^{n}-2u_{i}^{n}%
+u_{i-1}^{n}]+(1-\theta)\alpha_{i}^{n+1}[u_{i+1}^{n+1}-2u_{i}^{n+1}%
+u_{i-1}^{n+1}],
\end{array}
\label{SFDD-1}%
\end{equation}

where
\begin{equation}%
\begin{array}
[c]{rl}%
\alpha_{i}^{n}:= & \frac{\sigma(\Gamma_{i}^{n})^{2}\triangle t_{i}}%
{2(x_{i+1}-x_{i})(x_{i}-x_{i-1})},\\
\sigma(\Gamma_{i}^{n}):= & \left\{
\begin{array}
[c]{l}%
\underline{\sigma},\text{ \ \ if }\Gamma_{i}^{n}\geq0\\
\bar{\sigma},\text{ \ \ if }\Gamma_{i}^{n}<0
\end{array}
\right.  ,\\
\Gamma_{i}^{n}:= & \frac{u_{i+1}^{n}-2u_{i}^{n}+u_{i-1}^{n}}{(x_{i+1}%
-x_{i})(x_{i}-x_{i-1})}.
\end{array}
\end{equation}

In this paper, we consider the fully implicit schemes\ with $\theta=0,$ i.e.,%

\begin{equation}%
\begin{array}
[c]{cl}%
u_{i}^{n+1}-u_{i}^{n}= & \alpha_{i}^{n+1}[u_{i+1}^{n+1}-2u_{i}^{n+1}%
+u_{i-1}^{n+1}],
\end{array}
\label{SPDD-3}%
\end{equation}

The set of algebraic equation (\ref{sequa-3}) is nonlinear for the formula of
$\sigma(\Gamma_{i}^{n}).$ So we consider the discrete equation at each node as%

\[%
\begin{array}
[c]{c}%
\varphi_{i}^{n}:=u_{i}^{n}-u_{i}^{n+1}+\alpha_{i+1}^{n+1}[u_{i+1}^{n+1}%
-2u_{i}^{n+1}+u_{i-1}^{n+1}]
\end{array}
\]

For further analysis the Newton iteration, we rewrite the discrete equation
(\ref{SFDD-1}) in matrix form. Let%

\[%
\begin{array}
[c]{l}%
U^{n+1}=[u_{0}^{n+1},u_{1}^{n+1},\cdots,u_{m}^{n+1}]^{\prime},\text{ \ }%
U^{n}=[u_{0}^{n},u_{1}^{n},\cdots,u_{m}^{n}]^{\prime},\\
\lbrack M^{n}U^{n}]_{i}:=-\alpha_{i}^{n}[u_{i+1}^{n}-2u_{i}^{n}+u_{i-1}^{n}].
\end{array}
\]

\begin{equation}
\lbrack I+M^{n+1}]U^{n+1}=U^{n}, \label{SPDD-4}%
\end{equation}

\begin{theorem}
\label{w6} The fully implicit discretization (\ref{SPDD-4}) convergences to
the solution of the equation (\ref{sequa-3}), as $\triangle t$, $\triangle
x\rightarrow0.$
\end{theorem}

\renewcommand{\refname}


{\large References}

\begin{thebibliography}{99}                                                                                               %


\bibitem {Barles}G. Barles.(1997) \emph{Convergence of numerical schemes for
degenerate parabolic equations arising in finance,} In L. C. G. Rogers and
D.Talay (Eds), Numerical Methods in Finance, 1-21.

\bibitem {DHP11}L. Denis, M. Hu and S. Peng.(2011) \emph{Function spaces and
capacity related to a sublinear expectation: application to $G$-Brownian
motion paths,} Potential Anal, 139-161.

\bibitem {Pooey}D.M. Pooey, P.A. Forsyth and K.R. Vetzal.(2003)
\emph{Numerical convergence properties of option pricing PDEs with uncertain
volatility.} IMA. J. of Numerical Analysis 23, 241-267.

\bibitem {Peng2005}S. Peng.(2005) \emph{Nonlinear expectations and nonlinear
Markov chains,} Chin. Ann. Math. 26B(2), 159--184.

\bibitem {P07a}S. Peng.(2007) \emph{$G$-expectation, $G$-Brownian Motion and
Related Stochastic Calculus of It\^{o} type}, Stochastic analysis and
applications, Abel Symp., 2, Springer, Berlin. 541-567.

\bibitem {P08a}S. Peng.(2008) \emph{Multi-Dimensional $G$-Brownian Motion and
Related Stochastic Calculus under $G$-Expectation}, Stochastic Processes and
their Applications, 118(12), 2223-2253.

\bibitem {P09}S. Peng.(2009) \emph{Survey on normal distributions, central
limit theorem, Brownian motion and the related stochastic calculus under
sublinear expectations}, Science in China Series A: Mathematics, 52(7), 1391-1411.

\bibitem {Pang}J.S. Pang, L. Qi.(1993) \emph{Nonsmooth equations: Motivation
and algorithems.} SIAM Journal on Optimization 3, 443-465.

\bibitem {Qi}L. Qi, J. Sun.(1993) \emph{A nonsmooth version of Newton's
method.} Mathematical Programming 58, 353-367.

\bibitem {Sun}D. Sun, J. Han.(1997) \emph{Newton and quasi-Newton methods for
a class of nonsmooth equations and related problems.} SIAM Journal on
Numerical Analysis 17, 33-38.
\end{thebibliography}
\end{document}